\documentclass[graybox]{svmult}

\usepackage{mathptmx}       % selects Times Roman as basic font
\usepackage{helvet}         % selects Helvetica as sans-serif font
\usepackage{courier}        % selects Courier as typewriter font
\usepackage{type1cm}        % activate if the above 3 fonts are
\usepackage{makeidx}         % allows index generation
\usepackage{graphicx}        % standard LaTeX graphics tool
\usepackage{multicol}        % used for the two-column index
\usepackage[bottom]{footmisc}% places footnotes at page bottom

\makeindex             % used for the subject index
\usepackage[psamsfonts]{amsfonts}
\usepackage{amsmath,amssymb}
\usepackage{bbm} % used in \newcommand{\1}{\mathbbm{1}}

\newcommand{\Ebold} {{\mathbb E}}

\newcommand{\Nbold} {{\mathbb N}}

\newcommand{\Rbold} {{\mathbb R}}

\newcommand{\Zbold} {{\mathbb Z}}

\newcommand{\Scal}   {\mathcal{S}}

\newcommand{\Wcal}   {\mathcal{W}}

\newcommand{\R}{\Rbold}

\newcommand{\N}{\Nbold}

\newcommand{\1}{\mathbbm{1}}

\newcommand{\Zd}{\Zbold^{d}}
\newcommand{\nnb}   {\nonumber \\}
\begin{document}

\title*{The strong interaction limit of
        \\
        continuous-time weakly
        self-avoiding walk
         }
\author{David C.\ Brydges,
Antoine Dahlqvist,
and
Gordon Slade
}
\institute{
David C.\ Brydges \at
Department of Mathematics, University of British Columbia,
Vancouver, BC, Canada V6T 1Z2, \email{db5d@math.ubc.ca}
\and
Antoine Dahlqvist \at DMA--Ecole Normale Sup\'erieure, 45 rue d'Ulm, 75230 Paris
Cedex 5, France, \email{dahlqvis@clipper.ens.fr}
\and
Gordon Slade \at
Department of Mathematics, University of British Columbia,
Vancouver, BC, Canada V6T 1Z2, \email{slade@math.ubc.ca}
}
\maketitle

\vspace{-10mm} April 18, 2011

\bigskip \noindent
\emph{Dedicated to Erwin Bolthausen and J\"urgen G\"artner on the
occasion of their 65th and 60th birthday celebration}

\bigskip

\abstract{The strong interaction limit of the discrete-time weakly self-avoiding
walk (or Domb--Joyce model) is trivially seen to be the usual strictly
self-avoiding walk.  For the continuous-time weakly self-avoiding
walk, the situation is more delicate, and is clarified in this paper.
The strong interaction limit in the continuous-time setting depends on
how the fugacity is scaled, and in one extreme leads to the strictly
self-avoiding walk, in another to simple random walk.  These two
extremes are interpolated by a new model of a self-repelling walk that
we call the ``quick step'' model.  We study the limit both for walks
taking a fixed number of steps, and for the two-point function.}

\section{Domb--Joyce model:  discrete time}

The discrete-time weakly self-avoiding walk, or Domb--Joyce
model \cite{DJ72}, is a useful adaptation of the strictly
self-avoiding walk that continues to be actively studied
\cite{BRR11}.
It is defined as follows.  For simplicity, we restrict attention
to the nearest-neighbour model on
$\Zd$, although a more general formulation is easy to obtain.

Let $d \ge 1$ and $n \ge 0$ be integers, and let
$\Wcal_n$ denote the set of nearest-neighbour walks in $\Zd$,
of length $n$, which start from the origin.  In other words,
$\Wcal_n$ consists of sequences $Y=(Y_0,Y_1,\ldots,Y_n)$ with
$Y_i\in \Zd$, $Y_0=0$, $|Y_{i+1}-Y_i|=1$ (Euclidean distance).
Let $\Scal_n$ denote the set of nearest-neighbour
self-avoiding walks in $\Wcal_n$; these are the walks with
$Y_i \ne Y_j$ for all $i \ne j$.
Let $c_n$ denote the
cardinality of $\Scal_n$.
For $Y \in \Wcal_n$ and $x\in \Zd$, let
$n_x=n_x(Y) = \sum_{i=0}^n \1_{Y_i=x}$
denote the number of visits to $x$ by $Y$.
The Domb--Joyce model is the measure on $\Wcal_n$
which assigns to a walk
$Y \in \Wcal_n$ the probability
\begin{equation}
    P_{g,n}^{\rm DJ}(Y)
    =
    \frac{1}{c_n^{\rm DJ}(g)}e^{-g\sum_{x\in\Zd} n_x(Y)(n_x(Y)-1)},
\end{equation}
where $g$ is a positive parameter and
\begin{equation}
\label{e:cnDJdef}
    c_n^{\rm DJ}(g) = \sum_{Y\in \Wcal_n} e^{-g\sum_{x\in\Zd} n_x(Y)(n_x(Y)-1)}.
\end{equation}
The Domb--Joyce model interpolates between simple random walk
and self-avoiding walk.  Indeed, the case $g=0$ corresponds
to simple random walk by definition, and also
\begin{equation}
\label{e:limeg}
    \lim_{g \to \infty} e^{-g\sum_{x\in\Zd} n_x(Y)(n_x(Y)-1)}
    =
    \1_{Y \in \Scal_n}
\end{equation}
and hence
\begin{equation}
\label{e:limPn}
    \lim_{g \to \infty}P_{g,n}^{\rm DJ}(Y)
    =
    \frac{1}{c_n} \1_{Y \in \Scal_n}.
\end{equation}
This shows that the strong interaction limit of the Domb--Joyce model
is the uniform measure on $\Scal_n$.
(For an analogous result for weakly self-avoiding
lattice trees, which is more subtle
than for self-avoiding walks, see \cite{BCHS99}.)

A standard subadditivity argument
(see, e.g., \cite[Lemma~1.2.2]{MS93}) implies that the limits
\begin{equation}
    \mu(g) = \lim_{n \to \infty}c_n^{\rm DJ}(g)^{1/n},
    \quad\quad
    \mu = \lim_{n \to \infty} c_n^{1/n}
\end{equation}
exist and obey $c_n^{\rm DJ}(g) \ge \mu(g)^n$ and $c_n \ge \mu^n$ for
all $n$.
The number of walks that take steps only in the positive coordinate
directions is $d^n$, and such walks are self-avoiding, so $c_n \ge d^n$,
Also, it follows from \eqref{e:cnDJdef} that if $0 \le g < g_0$ then
$(2d)^n \ge c_n^{\rm DG}(g) \ge  c_n^{\rm DG}(g_0) \ge c_n \ge d^n$,
and hence $2d \ge \mu(g) \ge \mu(g_0) \ge \mu \ge d$.
In particular, by monotonicity, $\lim_{g \to \infty}
\mu(g)$ exists in $[\mu,2d]$.
If we
take the limit $g\to\infty$ in the inequality
$c_n^{\rm DJ}(g) \ge \mu(g)^n \ge \mu^n$, we obtain $c_n \ge
(\lim_{g \to\infty}\mu(g))^n
\ge \mu^n$.  Taking $n^{\rm th}$ roots and then the limit $n\to\infty$
then gives
\begin{equation}
\label{e:limmug}
    \lim_{g \to\infty}\mu(g)=\mu.
\end{equation}

Let $\Wcal_n(x)$ denote the subset of $\Wcal_n$ consisting
of walks that end at $x\in\Zd$.
Let $\Scal_n(x) = \Scal_n \cap \Wcal_n(x)$,
and let $c_n(x)$ denote the cardinality of $\Scal_n(x)$.
Let
\begin{equation}\
    c_{n,g}^{\rm DJ}(x)
    = \sum_{Y\in \Wcal_n(x)}
    e^{-g\sum_{\tilde{x}\in\Zd} n_{\tilde{x}}(Y)(n_{\tilde{x}}(Y)-1)}.
\end{equation}
Let $z \ge 0$.
The two-point functions of the Domb--Joyce and self-avoiding
walk models are defined as follows:
\begin{align}
    G_{g,z}^{\rm DJ}(x)
    &=
    \sum_{n=0}^\infty
    c_{n,g}^{\rm DJ}(x) z^n,
    \quad\quad
    G_z (x)
    =
    \sum_{n=0}^\infty
    c_n(x)  z^n.
    \label{e:Gdef}
\end{align}
These series converge for $z < \mu(g)^{-1}$ and $z < \mu^{-1}$
respectively.  Presumably they converge also for
$z = \mu(g)^{-1}$ and $z = \mu^{-1}$ but this is a delicate
question that is unproven except in high dimensions
(in fact, the decay of the two-point function with $z=\mu^{-1}$
is known in some cases
\cite{BS11,Hara08,HHS03}).
The following proposition shows that the strong interaction limit of
$G_{g,z}^{\rm DJ}(x)$ is $G_{z}(x)$.

\begin{proposition}
\label{prop:GDJ}
For $z \in [0, \mu^{-1})$ and $x \in \Zd$,
\begin{equation}
    \lim_{g \to \infty}  G_{g,z}^{\rm DJ}(x)
    = G_z(x).
\end{equation}
\end{proposition}

\begin{proof}
Fix $z \in [0, \mu^{-1})$. By \eqref{e:limmug}, if $g_0$ is sufficiently
large then $z < \mu(g_0)^{-1}$.  Thus, since $c_n^{\rm DJ}(g)$ is
nonincreasing in $g$,
there are $r<1$ and $C>0$ such that
$c_n^{\rm DJ}(g) z^n \le c_n^{\rm DJ}(g_0) z^n \le Cr^n$
for all $n$, uniformly in $g\ge g_0$.  Thus, for all $g\ge g_0$,
\begin{equation}
    G_{g,z}^{\rm DJ}(x)
    \le
    \sum_{x \in \Zd} G_{g,z}^{\rm DJ}(x)
    =
    \sum_{n=0}^\infty
    c_n^{\rm DJ}(g) z^n \le \frac{C}{1-r} < \infty.
\end{equation}
By \eqref{e:limeg}, $\lim_{g \to \infty} c_{n,g}^{\rm DJ}(x) = c_n(x)$,
and the desired result then follows by dominated
convergence.
\qed
\end{proof}

\section{The continuous-time weakly self-avoiding walk}

Our goal is to study the analogues of \eqref{e:limPn}
and Proposition~\ref{prop:GDJ} for the
continuous-time weakly self-avoiding walk.
The continuous-time model
is a lattice version of the Edwards model \cite{Edwa65}.
It has been useful in particular due to its
representation in terms of functional integrals \cite{BIS09}
that have been employed in renormalisation group analyses.

\subsection{Fixed-length walks}

We first consider the case of fixed-length walks, in which a
fixed number $n$ of steps is taken by the walk.
We will find that the strong interaction limit depends on how an
auxiliary parameter $\rho$ is scaled, where $e^\rho$ plays the role of a
fugacity. The scaling is parametrized by $a \in [-\infty
,\infty]$.  The case $a = \infty$ leads to the strictly self-avoiding
walk, the case $a = -\infty$ leads to simple random walk, and the
interpolating cases, $a \in (-\infty ,\infty)$, define a new
model of a self-repelling walk that we call the ``quick step'' model.

Let $X$ denote the continuous-time Markov process with state space
$\Zd$, in which uniformly random
nearest-neighbour steps are taken after independent
${\rm Exp}(1)$ holding times.  Let $\Ebold$ denote expectation for
this process started at $0$.  We distinguish between the
continuous-time walk $X$ and the sequence of sites visited during its
first $n$ steps, which we typically denote by $Y \in \Wcal_n$.
Conditioning on the first $n$ steps of $X$ to be $Y$ is denoted by
$\Ebold (\cdot \mid Y)$.

For fixed-length walks, the continuous-time weakly self-avoiding
walk is the measure $Q_{g,\rho,n}$ on $\Wcal_n$ defined as follows.
Here $\rho$ is a real parameter at our disposal,
which we allow to depend on $g>0$.
Let $T_n$ denote the time of the $(n+1)^{\rm st}$ jump of $X$,
and let $L_{x,n}(X) = \int_0^{T_n}\1_{X(s)=x}ds$ denote the local
time at $x$ up to time $T_n$.
By definition, $\sum_{x\in\Zd} L_{x,n}=T_n$.
For $Y\in\Wcal_n$, let
\begin{equation}
    Q_{g,\rho,n}(Y)
    =
    \frac{1}{Z_n(g,\rho)}
    \Ebold \left( e^{-g\sum_x L_{x,n}^2 + \rho \sum_x L_{x,n}}
    \mid Y
    \right),
\end{equation}
where
\begin{equation}
    Z_n(g,\rho) = \sum_{Y \in \Wcal_n}
    \Ebold \left( e^{-g\sum_x L_{x,n}^2 + \rho \sum_x L_{x,n}}
    \mid Y
    \right).
\end{equation}
For $a \in \R$ and $m \in \N$, let
\begin{equation}
\label{e:Indef}
    I_m(a) =
    \int_{-a}^\infty
    \frac{(a +u )^{m -1}}{(m -1)!}  e^{-u ^2}\, du.
\end{equation}

\begin{proposition}
\label{prop:fl} Let $\alpha = \alpha (g,\rho)= \frac 12
g^{-1/2}(\rho-1)$, and let $\rho=\rho(g)$ be chosen in such a way that
$a=\lim_{g \to \infty}\alpha(g,\rho (g))$ exists in
$[-\infty,\infty]$.  Let $n \ge 1$ and $Y \in \Wcal_n$.  Then
\begin{equation}
    \lim_{g \to \infty} Q_{g,\rho(g),n}(Y)
    =
    \begin{cases}
    \frac{1}{Z_a}\prod_{x \in Y}
    e^{a^2}
    I_{n_x(Y)}(a)
    &
    \text{if $a \in (-\infty,\infty)$,}
    \\
    \frac{1}{c_n} \1_{Y \in \Scal_n}
    &
    \text{if $a= \infty$,}
    \\
    \frac{1}{(2d)^n}
    &
    \text{if $a= -\infty$,}
    \end{cases}
\end{equation}
where $Z_a$ is a normalisation constant, and the product over $x$ is
over the distinct vertices visited by $Y$.
\end{proposition}

\begin{proof}
As before, we write $n_x = n_x(Y)$
for the number of times that $x$ is visited  by $Y$.
Thus $\sum_x n_x=n+1$ is the number of vertices visited by $Y$ (with
multiplicity).
Since the sum of $m$ independent ${\rm Exp}(1)$ random variables
has a Gamma$(m,1)$ distribution, we have
\begin{equation}
    \label{e:EY-1}
    \Ebold \left( e^{-g\sum_x L_{x,n}^2 + \rho \sum_x L_{x,n}}
    \mid Y
    \right)
    =
    \prod_{x \in Y}
    \int_0^\infty \frac{s_x^{n_x-1}}{(n_x-1)!} e^{-s_x}
    e^{-gs_x^2 + \rho s_x}\,ds_x,
\end{equation}
where the product is over the \emph{distinct} vertices visited by $Y$.
We make the changes of variables $t_x=g^{1/2}s_x$ and then
$u_x= t_x -\alpha$.
After completing the square, this leads to
\begin{align}
\label{e:EY0}
    \Ebold \left( e^{-g\sum_x L_{x,n}^2 + \rho \sum_x L_{x,n}}
    \mid Y
    \right)
    &=
    g^{-(n+1)/2}
    \prod_{x \in Y}
    e^{\alpha^2}
    I_{n_x}(\alpha)
    .
\end{align}

\medskip \noindent \emph{Case~$a\in(-\infty,\infty)$:
the quick step model}.
Suppose that $\alpha \to a \in (-\infty,\infty)$ as $g \to \infty$.
In this case, by the continuity of $I_m(a)$ in $a$,
\begin{equation}
    \Ebold \left( e^{-g\sum_x L_{x,n}^2 + \rho(g) \sum_x L_{x,n}}
    \mid Y
    \right)
    \sim
    g^{-(n+1)/2}
    \prod_{x \in Y}
    e^{a^2}
    I_{n_x}(a) ,
\end{equation}
and thus
\begin{equation}
    \label{e:limQ}
    \lim_{g \to \infty} Q_{g,\rho(g),n}(Y)
    =
    \frac{1}{Z_a}\prod_{x \in Y}
    e^{a^2}
    I_{n_x(Y)}(a)
    \quad\quad
    (\alpha \to a \in (-\infty,\infty)).
\end{equation}

\medskip \noindent \emph{Case~$a=\infty$: limit is uniform on} $\Scal_n$.
Suppose that $\alpha \to \infty$ as $g \to \infty$.
In this case, since $\alpha$ is nonzero we can use \eqref{e:EY0}
to write
\begin{align}
\label{e:EY}
&
    \Ebold \left( e^{-g\sum_x L_{x,n}^2 + \rho \sum_x L_{x,n}}
    \mid Y
    \right)
\nnb
    &=
    (g^{-1/2}e^{\alpha^2})^{n+1} (\alpha e^{-\alpha^2})^{n+1-|Y|}
    \prod_{x \in Y}
    \int_{-\alpha}^\infty
    \frac{(1+u_x/\alpha)^{n_x-1}}{(n_x-1)!}  e^{-u_x^2}\, du_x
    ,
\end{align}
where $|Y|$ denotes the number of distinct vertices visited by $Y$.
Since the factor $(\alpha e^{-\alpha^2})^{n+1-|Y|}$ goes to zero
unless $Y$ is self-avoiding, in which case the factor is equal to $1$
and $n_x=1$ for the vertices visited by $Y$, and since also
\begin{equation}
\label{e:J1lim}
    \lim_{\alpha \to \infty}\int_{-\alpha}^\infty
    e^{-u_x^2}\, du_x
    = \sqrt{\pi},
\end{equation}
this gives
\begin{align}
\label{e:EY2}
    \Ebold \left( e^{-g\sum_x L_{x,n}^2 + \rho(g) \sum_x L_{x,n}}
    \mid Y
    \right)
    &
    \sim
    (g^{-1/2}e^{\alpha^2}\sqrt{\pi})^{n+1} \1_{Y \in \Scal_n}.
\end{align}
When we take the normalisation into account we find that
\begin{equation}
    \lim_{g \to \infty} Q_{g,\rho(g),n}(Y)
    =
    \frac{1}{c_n} \1_{Y \in \Scal_n}
    \quad\quad
    (\alpha \to \infty).
\end{equation}

\medskip \noindent \emph{Case~$a=-\infty$: limit is uniform on} $\Wcal_n$.
Suppose that $\alpha \to -\infty$ as $g \to \infty$.
We will show that, for $m \ge 1$,
\begin{equation}
\label{e:aminus}
    e^{\alpha^{2}}I_{m}(\alpha)\sim (-2\alpha)^{-m}
    \quad \text{as $\alpha \to -\infty$}.
\end{equation}
With \eqref{e:EY0}, this claim implies that
\begin{equation}
    \Ebold \left( e^{-g\sum_x L_{x,n}^2 + \rho(g) \sum_x L_{x,n}}
    \mid Y
    \right)
    \sim
    g^{-(n+1)/2}
    \prod_{x \in Y}
    (-2\alpha)^{-n_x}
    = (-2\alpha g^{-1/2})^{n+1}
    .
\end{equation}
Since the right-hand side is independent of $Y$, this proves that the
limiting measure is uniform on $\Wcal_n$, as required.
Finally, to prove \eqref{e:aminus}, we set $b=-\alpha$ and obtain
\begin{align}
&
    (2b)^{m} e^{b^{2}} I_{m} (-b)
\nnb
    &=
    (2b)^{m} e^{b^{2}}
    \int_{b}^{\infty}\frac{(-b+u)^{m-1}}{(m-1)!} e^{-u^{2}}\,du \,
    =
    \int_{0}^{\infty}\frac{u^{m-1}}{(m-1)!} e^{- (u/ (2b))^{2} - u}\,du
    .
\end{align}
By dominated convergence, as $b \to \infty$, the integral on the
right-hand side approaches $1$ because it becomes the integral over
the $\Gamma (m,1)$ probability density function.
\qed\end{proof}

Proposition~\ref{prop:fl} shows that
the case $\alpha \to \infty$ leads to the uniform
measure on self-avoiding walks, whereas $\alpha \to -\infty$ leads
to simple random walk.  These two extremes are interpolated by
the quick step walk, for $\alpha \to a \in (-\infty,\infty)$
(e.g., $a=0$ if $|\rho| = o(g^{1/2})$
or $a=c$ if $\rho \sim 2c g^{1/2}$).
The name ``quick step walk'' is intended to reflect that idea
that the large $g$ limit of the continuous-time walk should be dominated
by quickly moving continuous-time walks.  In fact,
when $\rho = 2ag^{1/2}$, by completing the square
the weight $e^{-\sum_x (gL_{x,n}^2 - \rho L_{x,n})}$ can be
rewritten as $e^{\sum_x [-(g^{1/2}L_{x,n}-a)^2+a^2]}$.
Thus walks with smaller $L_{x,n}$ receive larger
weight, and this effect grows in importance as $g \to \infty$.

The particular case of Proposition~\ref{prop:fl} for the choice
\begin{equation}
\label{e:BFFrho}
    \rho(g) = \left(2g \log (g/\pi) \right)^{1/2},
\end{equation}
which corresponds to $a=\infty$, was proved previously in \cite{BFF84}.

For the case $a = 0$, evaluation of $I_{n_{x} (Y)}(0)$ in \eqref{e:limQ}
gives
\begin{equation}
    \lim_{g \to \infty} Q_{g,\rho(g),n}(Y)
    =
    \frac{1}{Z_0}
    \prod_{x \in Y}
    \frac{\Gamma(n_x(Y)/2)}{2\Gamma(n_x(Y))}
    \quad\quad
    (\alpha \to 0).
\end{equation}
Large values of $n_x$ are penalised under this limiting probability, so
this is a model of a self-repelling walk.
It is an interesting question whether
the quick step walk is in the same universality class as the
self-avoiding walk, for $a \in (-\infty,\infty)$.  We do not have
an answer to this question.

\subsection{Two-point function}

Now we show that when $\rho$ is chosen carefully, depending on $g$,
the two-point function for the continuous-time weakly self-avoiding
walk converges, as $g\rightarrow \infty $, to the two-point
function of the strictly self-avoiding walk.  The two-point function
of the continuous-time weakly self-avoiding walk can be written in two
equivalent ways.  This is discussed in a self-contained manner in
\cite{BIS09}, and we summarise the situation as follows.

The version of the two-point function that we will work with is
written in terms of a modified Markov process
$X=X(t)$, whose definition depends on a choice of $\delta \in (0,1)$.
The state space is $\Zd \cup \{\partial\}$, where $\partial$ is an
absorbing state called the cemetery.  When $X$ arrives at state $x$
it waits for an ${\rm Exp}(1)$ holding time and then jumps to a
neighbour of $x$ with probability $(2d)^{-1}(1-\delta)$ and jumps
to the cemetery with probability $\delta$.  The holding times
are independent of each other and of the jumps.
The two-point function is
defined, for $x \in \Zd$, to be
\begin{equation}
\label{e:GCTdef1}
    G_{g,\rho}^{\rm CT}(x) = \frac{1}{\delta}
    \Ebold^{(\delta)}
    \left( e^{-g\sum_{v\in\Zd} L_v^2 + \rho \zeta} \1_{X(\zeta^-)=x}
    \right),
\end{equation}
where we leave implicit the dependence of $G^{\rm CT}$ on $\delta$,
where $\Ebold^{(\delta)}$ denotes expectation with respect to
the modified process, and where $\rho$ is any real number for which the
expectation is finite.
The random number of steps taken by $X$ before jumping to the cemetery
is denoted $\eta$, and the independent sequence of holding times
will be denoted $\sigma_0,\sigma_1,\ldots,\sigma_\eta$.

A special case of the conclusions of \cite[Section~3.2]{BIS09}
(there with $d_x=1$ and $\pi_{x,\partial}=\delta$ for all $x$,
and restricted to finite state space)
is the equivalent formula
\begin{equation}
\label{e:GCTdef2}
    G_{g,\rho}^{\rm CT}(x) =
    \int_0^\infty
    \Ebold \left( e^{-g \sum_{v\in\Zd}L_{v,T}^2}\1_{X(T)=x} \right)
    e^{(\rho-\delta) T} dT,
\end{equation}
where now $X$ is the original continuous-time Markov process $X$
without cemetery state, and $\Ebold$ denotes its expectation
when started from the origin of $\Zd$.
Here $L_{v,T}=\int_0^T\1_{X(s)=v}ds$ is the local time of $X$ at $v\in\Zd$
up to time $T$.
We will work with \eqref{e:GCTdef1} rather than \eqref{e:GCTdef2}.

As in Proposition~\ref{prop:fl}, we write
$\alpha  = \alpha (g,\rho) = \frac 12 g^{-1/2}(\rho-1)$.
Throughout this section, we mainly choose $\rho=\rho(g)$
in such a way
that
\begin{equation}
\label{e:nu}
    \lim_{g \to \infty}g^{-1/2}e^{\alpha^2 (g,\rho (g))} = p \in [0,\infty)
\end{equation}
For example, \eqref{e:nu} holds for $p>0$ when
$\rho(g) = 2 [g \log (p \sqrt{g})]^{1/2}$, which is a choice closely
related to that in \eqref{e:BFFrho}.
Note that $\lim_{g\to\infty}\rho (g) = \infty$ when $p>0$.  It is
natural to consider $\rho \to \infty$, because if $\rho$ is fixed to a
value such that $G_{g_0,\rho}^{\rm CT}(x) < \infty$ for some $g_0>0$,
then by dominated convergence $\lim_{g\to\infty}G_{g,\rho}^{\rm
CT}(x) =0$.  The conclusion of Proposition~\ref{prop:2ptfcn} shows
that this trivial behaviour persists even when $\rho(g) \to \infty$ in
such a way that $p=0$.

Given $p \in [0,\infty)$, let
\begin{equation}
\label{e:pz}
    z=(2d)^{-1}(1-\delta)p\sqrt{\pi} .
\end{equation}
The following proposition shows that, under the scaling \eqref{e:nu},
the strong interaction limit of the continuous-time weakly self-avoiding
walk two-point function is the two-point function of the strictly
self-avoiding walk defined in \eqref{e:Gdef}.

\begin{proposition}
\label{prop:2ptfcn}
Let $\delta\in (0,1)$, $z \in [0,\mu^{-1})$, and $x \in \Zd$.
Suppose that \eqref{e:nu} holds with the value of $p\in [0,\infty)$
specified by $z$ via \eqref{e:pz}.  Then
\begin{equation}
\label{e:2ptlim}
    \lim_{g \to\infty}G_{g,\rho(g)}^{\rm CT}(x) = p\sqrt{\pi} G_z(x).
\end{equation}
\end{proposition}

The proof of Proposition~\ref{prop:2ptfcn}
uses three lemmas, and we discuss these next.
For $m \in \N$ and $\alpha >0$, let
\begin{equation}
    J_m(\alpha)
    =
    \int_{-\alpha}^\infty
    \frac{(1+u/\alpha)^{m-1}}{(m-1)!}  e^{-u^2} \,du.
\end{equation}

\begin{lemma}
\label{lem:J}
Given any $\epsilon>0$ there exists $A_0>0$ such that
for all $\alpha \ge A \ge A_0$ and $m \ge 1$,
\begin{align}
\label{e:Jep}
    J_{m}(\alpha)
    & \le
    (1+\epsilon) J_{m}(A).
\end{align}
\end{lemma}

\begin{proof}
For $m \ge 2$, $J_m(\alpha)$ is a non-increasing function of
$\alpha\in (0,\infty)$
because
\begin{align}
    \frac{dJ_m(\alpha)}{d\alpha}
    & =
    - \frac{1}{(m-2)!}
    \int_{-\alpha}^\infty \frac{u}{\alpha^2}(1+u/\alpha)^{m-2}e^{-u^2}du
    \nnb & =
    - \frac{1}{(m-2)!}
    \bigg[\int_{\alpha}^\infty \frac{u}{\alpha^2}(1+u/\alpha)^{m-2}e^{-u^2}du
\nnb
    &\quad
    +\int_{0}^\alpha \frac{u}{\alpha^2}
    [(1+u/\alpha)^{m-2} - (1-u/\alpha)^{m-2}] e^{-u^2}du \bigg]
    \nnb &
    \le 0
\end{align}
(note that in the first line the contribution from differentiating the
limit of integration vanishes),
and thus \eqref{e:Jep} holds even with $\epsilon = 0$.
For the remaining case $m=1$, since $J_1$ is increasing and
$\lim_{\alpha \to \infty} J_1(\alpha) = \sqrt{\pi}$ (see \eqref{e:J1lim}),
given any $\epsilon >0$ there exists $A_0>0$ such that if $\alpha \ge
A \ge A_0$
then $1 \le J_1(\alpha)/J_1(A) \le 1+\epsilon$.
\qed\end{proof}

%Let $\chi^{\rm CT}(g,\rho) = \sum_{x \in \Zd} G_{g,\rho}^{\rm CT}(x)$.
%Then, using \eqref{e:GCTdef1}, we can write
%$\chi^{\rm CT}(g,\rho) = \sum_{n=0}^\infty w_n(g,\rho)$, where
Recall that $\eta$ is the random number of steps taken by $X$
before jumping to the cemetery state.
For $x \in \Zd$, let
\begin{align}
    \label{e:wndef-1}
    w_n(g,\rho;x)
    &=
    \frac{1}{\delta}
    \Ebold^{(\delta)}
    [e^{-g\sum_v L_v^2 + \rho \zeta}  \1_{X(\zeta^-) =x} \1_{\eta =n}],
    \\
    w_n(g,\rho)
    &=
    \frac{1}{\delta}
    \Ebold^{(\delta)}
    [e^{-g\sum_v L_v^2 + \rho \zeta}  \1_{\eta =n}].
\end{align}
Let $w_n(g;x) = w_n(g,\rho(g);x)$
and $w_n(g) = w_n(g,\rho(g))$ with $\rho(g)$ chosen according to
\eqref{e:nu}.

\begin{lemma}
\label{lem:Jpi}
Suppose that \eqref{e:nu} holds with $p>0$, and let $z$ be given
by \eqref{e:pz}.
Then for $n \ge 0$ and $x\in \Zd$,
\begin{equation}
    \lim_{g \to \infty}w_n(g;x)
    =
    p\sqrt{\pi} c_n(x)z^n.
\end{equation}
\end{lemma}

\begin{proof}
Given that $\eta =n$, let $Y\in\Wcal_n(x)$ denote the sequence of
jumps made by $X$ before landing in the cemetery, and let $|Y|$ denote
the cardinality of the range of $Y$.  By conditioning on $Y$ and using
\eqref{e:EY}, we see that, as $g \rightarrow \infty$,
\begin{align}
    w_n(g ;x)
    &
    =
    [(2d)^{-1}(1-\delta)]^n(g^{-1/2}e^{\alpha^2})^{n+1}
    \sum_{Y \in \Wcal_n(x)}
    (\alpha e^{-\alpha^2})^{n+1-|Y|}
    \prod_{v \in Y}
    J_{n_v}(\alpha)
    \nnb
    &
    \sim
    [(2d)^{-1}(1-\delta)]^n p^{n+1}
    \sum_{Y \in \Wcal_n(x)}
    (\alpha e^{-\alpha^2})^{n+1-|Y|}
    \prod_{v \in Y}
    J_{n_v}(\alpha)
    ,
\label{e:wndef}
\end{align}
where the product is over the distinct vertices visited by $Y$
and $|Y|$ denotes the number of such vertices.
It suffices to show that, for
any $Y \in \Wcal_n(x)$,
\begin{equation}
    \lim_{g \to \infty} (\alpha e^{-\alpha^2})^{n+1-|Y|}
    \prod_{v \in Y} J_{n_v}(\alpha)
    =
    \1_{Y \in \Scal_n} \pi^{(n+1)/2}.
\end{equation}
Since $p>0$, we have $\alpha \to \infty$, and so
$\alpha e^{-\alpha^2} \to 0$.
Therefore, the above limit is zero unless
$n+1=|Y|$, which corresponds to $Y \in \Scal_n$;
the product over $v$ remains bounded as $\alpha \to \infty$
and poses no difficulty.
Since $J_1(\alpha) \to \sqrt{\pi}$
as in \eqref{e:J1lim}, the result follows.
\qed\end{proof}

\begin{lemma}
Suppose that \eqref{e:nu} holds with $p \in (0,\infty)$,
and let $z$ be specified by \eqref{e:pz}.
\label{lem:wn}
Let
\begin{equation}
    \mu(g,\rho ) = \limsup_{n \to \infty} w_n(g,\rho )^{1/n}.
\end{equation}
Then
\begin{equation}
    \limsup_{g \to \infty} \mu(g,\rho (g) ) \le z\mu.
\end{equation}
\end{lemma}

\begin{proof}
Let $L_{x,[i,j]} = \sum_{k=i}^j \sigma_k \1_{Y_k=x}$, where the
$\sigma_k$ are the exponential holding times. Let
$\Ebold^{(\delta)}_{y}$ denote the expectation for the process started
in state $y$ instead of state $0$.
For integers $n\ge 1$ and $m\ge 1$,
an elementary argument using the
strong Markov property leads to

\begin{align}
    w_{n+m}(g,\rho ) & \le
    \frac{1}{\delta}
    \Ebold^{(\delta)}[e^{-g\sum_x L_{x,[0,n]}^2 + \rho \sum_x L_{x,[0,n]}}
    e^{-g\sum_x L_{x,[n+1,n+m]}^2 + \rho \sum_x L_{x,[n+1,n+m]}}
    \1_{\eta =n+m}]
    \nnb & =
    \sum_y
    \Ebold^{(\delta)}[e^{-g\sum_x L_{x,[0,n]}^2 + \rho \sum_x L_{x,[0,n]}}\1_{Y_{n+1}=y}]
    \frac{1}{\delta}
    \Ebold^{(\delta)}_y[e^{-g\sum_x L_{x}^2 + \rho \sum_x L_{x}}
    \1_{\eta =m-1}]
    \nnb & =
    \frac{1-\delta}{\delta}
    \Ebold^{(\delta)}
    [e^{-g\sum_x L_{x,[0,n]}^2 + \rho \sum_x L_{x,[0,n]}}\1_{\eta =n}]
    w_{m-1}(g,\rho)
    \nnb & \le
    w_n (g,\rho) w_{m-1}(g,\rho).
\end{align}
It is straightforward to adapt the proof of \cite[Lemma~1.2.2]{MS93}
to obtain from this approximate subadditivity the equality
\begin{equation}
    \mu(g,\rho)
    =
    \inf_{n \ge 1}w_n(g,\rho)^{1/(n+1)}.
\end{equation}
Then we have
\begin{equation}
\label{e:wmu}
    w_n(g,\rho)^{1/(n+1)} \ge \mu(g,\rho) .
\end{equation}

We let $g \to \infty$ in the above inequality,
with $\rho(g)$ chosen as in \eqref{e:nu};
note that $\alpha \to \infty$
since $p>0$.
%Given that $\eta =n$, let $Y\in\Wcal_n(x)$ denote the sequence of
%jumps made by $X$ before landing
%in the cemetery, and let $|Y|$ denote the cardinality of the range of $Y$.
%By conditioning on $Y$ and using \eqref{e:EY}, we see that
%\eq
%\label{e:wndef}
%    w_n(g,\rho)
%    =
%    [(2d)^{-1}(1-\delta)]^n(g^{-1/2}e^{\alpha^2})^{n+1}
%    \sum_{Y \in \Wcal_n}
%    (\alpha e^{-\alpha^2})^{n+1-|Y|}
%    \prod_{v \in Y}
%    J_{n_v}(\alpha).
%\en
%% For $p=0$, \eqref{e:wndef} implies that $\lim_{g \to \infty} w_n(g)
%% =0$.
%Since
%$\alpha \to \infty$, we can apply
By Lemma~\ref{lem:Jpi},
for $n \ge 0$,
\begin{equation}
\label{e:wncn}
    \lim_{g \to \infty} w_n(g)
    = p\sqrt{\pi} c_n z^n .
\end{equation}
By \eqref{e:wmu}, this gives
\begin{equation}
    (p\sqrt{\pi}c_n)^{1/(n+1)}z^{n/(n+1)}
    \ge
    \limsup_{g \to \infty} \mu(g,\rho(g)) .
\end{equation}
Now we take $n \to \infty$ to get
\begin{equation}
    \mu z \ge \limsup_{g \to \infty} \mu(g,\rho(g)),
\end{equation}
as required.
\qed\end{proof}

\medskip \noindent
\emph{Proof of Proposition~\ref{prop:2ptfcn}.}
We consider separately the cases $p>0$ and $p=0$.

\smallskip \noindent
\emph{Case $p>0$.}  We write $\rho = \rho(g)$.
By \eqref{e:GCTdef1}, and by \eqref{e:wndef-1} with $\rho=\rho(g)$,
\begin{equation}
    \label{e:eta-decomposition}
    G_{g,\rho}^{\rm CT}(x) =
    \sum_{n=0}^\infty
    w_n(g ;x)
    .
\end{equation}
By Lemma~\ref{lem:Jpi}, the result of taking the limit $g \to \infty$
under the summation gives the desired result
\begin{equation}
    p\sqrt{\pi} \sum_{n=0}^\infty c_n(x) z^n ,
\end{equation}
and it suffices to justify the interchange of limit and summation.
For this, we will use dominated convergence.  Since
$w_{n} (g;x) \le w_{n} (g)$, it suffices to find a
$g_{0}>0$ and a summable sequence $B_n$ such that, for $g
\ge g_{0}$ and $n \in \Nbold_{0}$,
\begin{equation}
    w_n(g;x)
    \le B_n.
\end{equation}
This will follow if we show the stronger statement that for large $g$
\begin{equation}
    w_n(g) \le B_n.
\end{equation}

Since $z \mu < 1$, there exists $\epsilon >0$ such that $c=
(1+\epsilon)^{2} (\mu z+\epsilon)<1$.  Since $g^{-1/2}e^{\alpha^2}
\rightarrow p>0$, there is a (large) $g_0$ such that if $g \ge g_0$
then $g^{-1/2}e^{\alpha^2} \le g_0^{-1/2}e^{\alpha_0^2}(1+\epsilon)$,
where $\alpha_0$ is the value of $\alpha$ corresponding to $g=g_0$;
also $\alpha e^{-\alpha^2} \le \alpha_0 e^{-\alpha_0^2}$.  Therefore,
by \eqref{e:wndef}, and by Lemma~\ref{lem:J} (increasing $g_0$ if
necessary),
\begin{align}
    w_n(g) &= [(2d)^{-1}(1-\delta)]^n (g^{-1/2}e^{\alpha^2})^{n+1}
    \sum_{Y \in \Wcal_n }
    (\alpha e^{-\alpha^2})^{n+1-|Y|}
    \prod_{v \in Y}
    J_{n_v}(\alpha)
    \nnb &  \le
    [(2d)^{-1}(1-\delta)]^n (g_0^{-1/2}e^{\alpha_0^2}(1+\epsilon)^2)^{n+1}
    \sum_{Y \in \Wcal_n }
    (\alpha_0 e^{-\alpha_0^2})^{n+1-|Y|}
    \prod_{v \in Y}
    J_{n_v}(\alpha_0)
    \nnb & =
    (1+\epsilon)^{2(n+1)} w_n(g_0).
\end{align}
We set $B_n = (1+\epsilon)^{2(n+1)}w_n(g_0)$.  Then
\begin{equation}
    \limsup_{n \to \infty} B_n^{1/n} = (1+\epsilon)^2\mu(g_0,\rho(g_0))
    \le (1+\epsilon)^2 (z\mu + \epsilon) < 1,
\end{equation}
by taking $g_0$ larger if necessary and applying Lemma~\ref{lem:wn}.
Therefore $\sum_n B_n$ converges, and the proof is complete for
the case $p>0$.

\smallskip \noindent
\emph{Case $p=0$.}  We will prove that
\begin{equation}
    \label{e:pzero-goal}
    \lim_{g\rightarrow \infty}
    \sum_{n = 0}^\infty w_n(g )
    =
    0
    .
\end{equation}
By \eqref{e:eta-decomposition}, this is more than sufficient.
We again write $\rho=\rho(g)$.
By conditioning on $Y$ and using
\eqref{e:EY0}, for $n \ge 0$ we have
\begin{align}
\label{e:wid}
    w_n(g)
    &=
    [(2d)^{-1}(1-\delta)]^n
    \sum_{Y \in \Wcal_n}
    \prod_{x \in Y}
    g^{-n_x/2} e^{\alpha^2}
    I_{n_x}(\alpha)
    .
\end{align}
%where
%\begin{equation}
%    K_{n} (g)
%    =
%    \int_0^\infty ds \frac{s^{n-1}}{(n-1)!} e^{-s} e^{-gs^2}
%    e^{\rho s}
%    .
%\end{equation}
The change of variables $s=a+u$ in \eqref{e:Indef} gives, for $m \ge 1$,
\begin{align}
    e^{\alpha^2} I_{m} (\alpha)
    &=
    e^{\alpha^{2}}
    \int_0^\infty \frac{s^{m-1}}{(m-1)!} e^{- (s-\alpha)^2}\,ds
    \nnb
    &\le
    e^{\alpha^{2}}
    \int_0^\infty \frac{s^{m-1}}{(m-1)!} e^{-s}
    \left(\sup_{s\in \Rbold}e^{s - (s-\alpha)^2}\right) ds
    =
    e^{\alpha^{2}+\alpha + 1/4}
    .
\label{e:Kbd}
\end{align}
Let $\epsilon >0$. Since $g^{-1/2}e^{\alpha^{2}}
\rightarrow p=0$, we can find $g (\epsilon)$ such that for $g\ge g
(\epsilon)$ and $m\ge 2$,
\begin{equation}
\label{e:gep}
    g^{-1/2}e^{\alpha^{2}}\sqrt{\pi}
    \le
    \epsilon,
    \quad \quad
    g^{-m/2}e^{\alpha^{2} + \alpha + 1/4}
    \le
    \epsilon^{m}
    .
\end{equation}
Henceforth we assume that $g \ge  g (\epsilon)$.
By \eqref{e:Kbd},
\begin{align}
    g^{-m/2} e^{\alpha^2}
    I_{m}(\alpha)
    &
    \le
    \epsilon^{m}
    \quad \text{for $m\ge 2$}
    .
\end{align}
For $m=1$, we obtain an upper bound by extending the range of the
integral in the first line of \eqref{e:Kbd}
to the entire real line, whereupon it
evaluates to $\sqrt{\pi}$.  Thus, by \eqref{e:gep},
$g^{-1/2}e^{\alpha^{2}} I_1(\alpha) \le \epsilon$.  By
\eqref{e:wid} and the fact that the number of walks in $\Wcal_n$ is
$(2d)^{n}$, for $n \ge 0$ we then have
\begin{align}
    w_n(g)
%    &=
%    [(2d)^{-1}(1-\delta)]^n
%    \sum_{Y \in \Wcal_n}
%    \prod_{v \in Y}
%    K_{n_v}(\alpha)
%    \nnb
    &\le
    [(2d)^{-1}(1-\delta)]^n
    \sum_{Y \in \Wcal_n}
    \prod_{v \in Y}
    \epsilon ^{n_v}
    =
    (1-\delta)^n
    \epsilon ^{n+1}
    .
\end{align}
(The case $n=0$ corresponds to $m=1$ because the number of visits
to state $0$ is $n_{0}=1$.)  Therefore $\limsup_{g \rightarrow
\infty}\sum_{n= 0}^\infty w_{n}(g) = O (\epsilon)$.  Since $\epsilon$
is arbitrary, this proves \eqref{e:pzero-goal}, and the proof is
complete.
\qed

\acknowledgement
The work of DB and GS was supported in part by NSERC of
Canada.

% \bibliography{../bibdef/bib}
% %\bibliography{references-my-papers}
% \bibliographystyle{plain}

\end{document}